\documentclass{amsart}
\usepackage[latin1]{inputenc}
\usepackage{amssymb,amsmath}

\newenvironment{dem}[1]%
    {\par\noindent{\scshape Proof:\ }#1}%
    {\mbox{}\hfill$\square$\par\bigskip\par}
\newtheorem{teo}{Theorem}
\newtheorem{defi}[teo]{Definition}
\newtheorem{prop}[teo]{Proposition}

\newtheorem{lema}[teo]{Lemma}

\newtheorem{prob}{Problem}
\newcommand{\To}{\longrightarrow}

\begin{document}

\title{Commutative rings with finite quotient fields}

\author{Antonio Avil\'es}
\thanks{Author supported by FPU grant of SEEU-MECD, Spain.}
\email{avileslo@um.es}
\address{Departamento de Matem\'aticas. Universidad de
Murcia, 30100 Murcia (Spain)}

\begin{abstract}We consider the class of all commutative reduced rings for which
there exists a finite subset $T\subset A$ such that all
projections on quotients by prime ideals of $A$ are surjective
when restricted to $T$. A complete structure theorem is given for this class of 
rings, and it is studied its relation with other finiteness conditions on 
the quotients of a ring over its prime ideals.\end{abstract}

\maketitle

\section*{Introduction}

Our aim is to study the structure of commutative rings that satisfy 
suitable finiteness conditions on its quotients by prime ideals. If $A$ is 
a commutative ring and $A/p$ represents the quotient ring of $A$ by a 
prime ideal $p$, we will be interested in the following conditions on 
$A$:\\
\begin{enumerate}
\item  All $A/p$ are finite.\\

\item The cardinal of all $A/p$ is bounded by some $n$.\\

\item There are $x_1,\ldots,x_n\in A$ such that $A/p = 
\{x_1+p,\ldots,x_n+p\}$ for any prime ideal $p$ of $A$.\\ 
\end{enumerate}

Clearly, $3\Rightarrow 2\Rightarrow 1$ and at the end of this article 
examples can be found showing that no converse is true. Also, observe that 
the three conditions are the same for $A$ as for the reduced ring 
$A/N(A)$, so we may restrict ourselves to reduced rings, that is, to commutative rings 
without nilpotent elements. One main result in this paper is that a complete structure 
theorem can be given for reduced rings satisfying condition (3). Namely, 
we associate to each such a ring $A$ a tuple $(K_1,B_1,\ldots,K_n,B_n)$, 
where the $K_i$'s are non isomorphic finite fields and the $B_i$'s are 
Boolean rings, in such a way that two of these rings are isomorphic if and 
only if the associated tuples are equal, up to isomorphism and order.\\

A reduced ring satisfying just condition (1) must be Von Neumann regular
(absolutely flat, in other terminology), as it is any reduced ring in 
which all prime ideals are maximal~\cite[Theorem 1.16]{Goodearl}. In this 
sense, our work continues that of  N. Popescu and C. Vraciu 
in~\cite{Vra1}, where the structure of commutative Von Neumann 
regular rings is studied. Here it is shown that imposing these kind of finiteness conditions,
strong structure results can be given. Let us note that C. Vraciu had 
already found in~\cite{Vra2}, in other direction, the following result concerning rings
satisfying condition (1):

\begin{teo}\label{Vraciu}
Let $B$ be a Boolean ring and $k:Spec(B)\To FFields$ a map from the Zariski spectrum $Spec(B)$ 
to the class of finite fields, such that for each 
$p\in Spec(B)$ there is a neighbourhood $U$ of $p$ such that $k(p)\subseteq k(p')$
for all $p'\in U$. Then, there is a commutative Von Neumann regular ring $A$ such that $B(A)=B$ 
and for each $p\in Spec(B)$, $A/pA$ is isomorphic to $k(p)$.
\end{teo}

On the other hand, the class of CFG-rings introduced in~\cite{Avi} is exactly the class of 
reduced rings satisfying condition (3), so the structure theorem for this 
class presented here is a precise measure of the degree of generality of the results 
in~\cite{Avi} in terms of CFG-rings. Also, a result in~\cite{Avi} suggests the addition
of the following condition to our list:\\

\renewcommand{\theenumii}{\roman{enumii}}
\renewcommand{\theenumi}{\arabic{enumi}.\arabic{subsection}}
\setcounter{subsection}{5}

\begin{enumerate}
\setcounter{enumi}{1}
 \item For any map $f:A\To A$ the following are equivalent
\begin{enumerate}
\item There exists a polynomial $F\in A[X]$ such that $f(x)=F(x)$ for all 
$x\in A$.
\item For all $x,y\in A$, $(f(x)-f(y))A\subseteq (x-y)A$.\\
\end{enumerate}
\end{enumerate}

\renewcommand{\theenumi}{\arabic{enumi}}

Lemma 3.3 in~\cite{Avi} is just the assertion that $(3)$ implies $(2.5)$ and Proposition 13 
below states that $(2.5)$ implies $(2)$. We give an example showing that 
the last converse is not true but we have been unable to decide about the 
first one.

\section*{Notations and terminology}

In the sequel, all rings will be supposed to be commutative with identity.
Letter $A$ will always represent a ring.
A reduced ring is a
ring without nilpotent elements.\\

A polynomial map $f:A^{n}\To A^{m}$ is a map for which all
components $f_{i}:A^{n}\To A$ are given by polynomials with
coefficients in $A$.\\

A Boolean ring is a ring $B$ such that $x^2=x$ for all $x\in B$.\\

For a Boolean ring $B$ and $a,b\in B$, the expression $a\leq b$ will mean 
$aB\subseteq bB$.\\

Let $A$ be a ring. We will denote by $B(A)$
the set of all idempotent elements of $A$. Recall that the
set $B(A)$ has a structure of Boolean ring with product inherited
from $A$ and with the sum $a\tilde{+}b=(a-b)^{2}$.\\

 A complete family of orthogonal
idempotents (shortly \textbf{c.f.o.i}) of $A$ is a family
$a_{1},\ldots,a_{n}$ of elements of $B(A)$ such that $a_{i}a_{j}
=0$ for $i\neq j$ and $\sum_{1}^{n}a_{i}=1$.\\

A convex combination of elements $x_{1},\ldots,x_{n}\in A$ is a linear combination 
$\sum_{1}^{n}a_{i}x_{i}$ such that $a_{1},\ldots,a_{n}$ is a c.f.o.i. of 
$A$. The set of all convex combinations of elements of a set $S\subseteq 
A$ will be denoted by $conv(S)\subseteq A$.\\

Regular ring will mean here
commutative Von Neumann regular ring (also called commutative absolutely flat rings), i.e. a (commutative) ring
for which any principal ideal is generated by an idempotent. If $A$ is a regular ring,
$e:A\longrightarrow B(A)$ will be the map that sends each $a\in A$
to the only idempotent $e(a)\in B(A)$ such that $aA=e(a)A$.\\

For a ring $A$, $Spec(A)$ will denote its Zariski spectrum,
that is, the topological space whose underlying set is the 
set of all prime ideals of $A$ and whose closed subsets are of the form 
$V(I) = \{ p : I\subseteq p\}$ being $I$ an ideal of $A$. If $A$ is a regular ring,
$B(A)$ can be 
identified with the ring of closed-open sets of $Spec(A)$ via the 
bijection $b\leftrightarrow O_{b} = \{p : b\not\in p\}$. In fact, for all 
$a\in A$, $O_{a}=O_{e(a)}$ is a closed and open set since its complement is 
$O_{1-e(a)}$. Also, recall that the correspondence $p\mapsto p\cap B(A)$ induces 
a homeomorphism $Spec(A)\To Spec(B(A))$, so that both spectra can be 
identified.

\section{Characterizations of CFG-rings}

\begin{teo}\label{caracterizacion de los anillos CFG}
Let $A$ be a reduced ring and let $x_1,\ldots,x_n\in A$. The following are equivalent:
\begin{enumerate}
\item $A = conv\{x_1,\ldots,x_n\}$.
\item $\prod_{i=1}^{n}(x-x_{i}) = 0$ for all $x\in A$.
\item $A/p = \{x_{1}+p,\ldots,x_{n}+p\}$ for every prime ideal $p$
of $A$.
\end{enumerate}
\end{teo}

\begin{dem}
    Suppose $x\in conv\{x_1,\ldots,x_n\}$, so 
    that $x=\sum_1^n a_i x_i$ and $\{a_i\}$ is a c.f.o.i. Then,
    $$\prod_{i=1}^n(x-x_i) = \prod_{i=1}^n\sum_{j=1}^n a_j(x_j-x_i)
    =\sum_{j=1}^n \prod_{i=1}^n a_j(x_j-x_i) = 0$$
    and this proves that $(1)$ implies $(2)$.\\
    
    Clearly $(2)$ implies $(3)$ 
    since the ideal $p$ is prime. If $(3)$ holds, then 
    $\prod_{i=1}^{n}(x-x_{i})$ belong to all prime ideals, so since $A$ is 
    reduced, this product is zero and $(2)$ holds. So $(2)$ and $(3)$ are 
    equivalent.\\
    
    Suppose now that $(2)$ holds.  Let 
    $x\in A$ and we show that $x\in conv\{x_1,\ldots,x_n\}$. We have the following 
    equality in $B(A)$:
    $$\prod_{i=1}^n e(x-x_i) = 0.$$
    Let $b_i=e(x-x_i)$ and $a_i = (1-b_i)\prod_{j=1}^{i-1}b_j$, then 
    $\{a_i\}_{i=1}^n$ is a c.f.o.i and
    $$\sum_{i=1}^n a_i e(x-x_i)=0$$
    so that $$e(x-\sum_{1}^{n}a_{i}x_{i}) = 
e(\sum_{1}^{n}a_{i}(x-x_{i})) = \sum_{1}^{n}a_{i}e(x-x_{i})=0.$$
and finally $x=\sum_1^n a_i x_i$.
\end{dem}

Observe that condition $(3)$ in Theorem~\ref{caracterizacion de los anillos CFG}
implies that each prime ideal is maximal, so $A$ is regular
by~\cite[Theorem 1.16]{Goodearl}. Following~\cite{Avi}, a CFG-ring is a 
regular ring satisfying condition (1) of this theorem for some 
$x_1,\ldots,x_n$. Equivalently, a CFG-ring is a reduced ring for which 
there are $x_1,\ldots,x_n$ satisfying any of the conditions of the 
theorem.

\section{Boolean envelopes of fields}

In this section we define a family of CFG-rings for which we
will prove in Theorem~\ref{estructura de CFG-anillos} that all
CFG-rings are finite products of rings of this type in a unique
way. Given a field $K$ and a Boolean ring $B$, it is defined a
$K$-algebra $K^{[B]}$ that will be a CFG-ring if $K$ is a finite
field. This algebras were firstly introduced by C. 
Vraciu~\cite{Vra2}.
Roughly speaking, $K^{[B]}$  will be the set of all
\emph{``formal convex combinations''} of elements of $K$ with
coefficients in $B$, where the sum and product are defined in
order to extend the operations in $K$ and to \emph{``commute with
convex combinations''}. More precisely:

\begin{defi}
Let $B$ be a Boolean ring and $K$ a field. The $K$-algebra of all 
continuous functions from $Spec(B)$ to $K$ ($K$ is equipped with the discrete 
topology) will be denoted by $K^{[B]}$.
\end{defi}

We make some elementary remarks:

\begin{enumerate}

\item We identify $K$ inside $K^{[B]}$, identifying each $x\in K$ with the 
corresponding constant map. In this way, $K^{[B]}$ is a $K$-algebra.\\

\item We identify $B$ inside $K^{[B]}$, identifying each $b\in B$ with the 
map $Spec(B)\To K$ that is constant equal 1 on $b$ and vanishes outside $b$. In fact, this 
gives all idempotent elements of $K^{[B]}$. Shortly, $B(K^{[B]})=B$.\\

\item The ring $K^{[B]}$ consists exactly of all functions $u:Spec(B)\To K$ with 
finite image, such that $u^{-1}(x)$ is a closed-open set of $Spec(B)$ for 
all $x\in K$. Making the identifications above, this means that each element of $K^{[B]}$ has an expression like 
$u=\sum_{1}^{n}a_{i}x_{i}$ where $x_{1},\ldots,x_{n}$ are elements 
of $K$ (the image of $u$) and $a_{1},\ldots,a_{n}$ is c.f.o.i of $B$ ($a_{i}=u^{-1}(\{x_i\})$). So $conv(K) = 
K^{[B]}$.\\

\item $K^{[B]}$ is a regular ring. Clearly, each element is the product of 
a unit (a $(K\setminus\{0\})$-valued continuous map) and an idempotent (a 
$\{0,1\}$-valued continuous map).\\

\item If $K$ is finite, then $K^{[B]}$ is a CFG-ring.\\ 
\end{enumerate}

\begin{prop}\label{unicidad de algebras regulares}
    Let $K$ be a field, $B$ a Boolean ring and $A$ a Von Neumann regular 
    commutative
    $K$-algebra such that $B(A)$ is isomorphic to $B$ and $A=conv(K)$. Then,
    $A$ and $K^{[B]}$ are isomorphic $K$-algebras.
\end{prop}

\begin{dem}
   We define a map $f:K^{[B]}\To A$ in the following way: if 
  $u=\sum_1^{n}a_{i}x_{i}\in K^{[B]}$ being $a_{1},\ldots,a_{n}$ a c.f.o.i 
  of $B$
  and $x_{i}\in K$, then $f(u) = \sum_1^{n}a_{i}x_{i}\in A$. 
  First of all, we must check that this is a correct definition, that is, 
  that $f(u)$ does not depend on the expression chosen for $u$ as convex 
  combination of elements of $K$. By elementary reduction arguments, it is enough to prove the following:
  \begin{itemize}
  \item[($\star$)] Let $a_{1},\ldots,a_{n}$ and $b_{1},\ldots,b_{n}$ be c.f.o.i. of 
  $B$, and let $x_{1},\ldots,x_{n}$ be $n$ different elements of $K$. If 
  $$\sum_{1}^{n}a_{i}x_{i} = \sum_{1}^{n}b_{i}x_{i}\in 
  K^{[B]}$$ then $$\sum_{1}^{n}a_{i}x_{i} = \sum_{1}^{n}b_{i}x_{i}\in 
  A.$$
  \end{itemize}
Under those hypotheses $u = \sum_{1}^{n}a_{i}x_{i}\in 
K^{[B]}$ is the map $u:Spec(B)\To K$ that takes the value $x_{i}$ on $a_{i}$. Therefore
$$\sum_{1}^{n}a_{i}x_{i} = \sum_{1}^{n}b_{i}x_{i}\in K^{[B]}\Rightarrow a_{i}=b_{i} \text{for all }i$$    
and assertion ($\star$) follows immediately.\\

We prove now that $f$ is a $K$-algebra 
homomorphism. We check for instance that $f$ commutes with addition (it is analogous for product).
We take $x=\sum_{i}a_{i}x_{i}$ and $y=\sum_{j}b_{j}y_{j}$ in $K^{[B]}$ 
being $a_{1},\ldots,a_{n}$ and $b_{1},\ldots,b_{m}$ c.f.o.i. of $B$ and 
$x_{i},y_{j}\in K$. Note that the $c_{ij} = a_{i}b_{j}$ constitute a 
c.f.o.i. of $B$, \begin{eqnarray*}
f(x+y) &=& f\left(\sum_{i}a_{i}x_{i}+\sum_{j}b_{j}y_{j}\right) = 
f\left(\sum_{i,j}c_{ij}x_{i}+\sum_{i,j}c_{ij}y_{j}\right) \\ &=& 
f\left(\sum_{i,j}c_{ij}(x_{i}+ y_{j})\right) = \sum_{i,j}c_{ij}(x_{i}+ 
y_{j}) \\ &=&
\sum_{i,j}c_{ij}x_{i}+\sum_{i,j}c_{ij}y_{j} = 
\sum_{i}a_{i}x_{i}+\sum_{j}b_{j}y_{j}\\ &=& f(x)+f(y)
\end{eqnarray*}
 On the other hand, $f$ is onto since $A=conv(K)$. We check that $f$ is one to one. Suppose 
 $f(u)=0$ and $u=\sum_{i}a_{i}x_{i}$ with the $a_{i}$'s a c.f.o.i of $B$ and 
 $x_{1}=0,x_{2}\neq 0,\ldots,x_{n}\neq 0$ elements of $K$. Then, since 
 $0=\sum_{i}a_{i}x_{i}\in A$, by multiplication by $a_{i}$, $a_{i}x_{i}=0$ 
 for all $i$. Since $x_{i}\in K\setminus\{0\}$ is a unit for all $i>1$, 
 $a_{i}=0$ for $i>1$. This implies that $u=0$.
\end{dem}

\section{Structure theorem for CFG-rings}

\begin{lema}\label{producto}
    Let $K$ be a finite field and $B_{1},B_{2}$ Boolean rings.
    Then, $K^{[B_{1}\times B_{2}]}$ and $K^{[B_{1}]}\times K^{[B_{2}]}$ 
    are isomorphic $K$-algebras.
\end{lema}

\begin{dem}
The ring $A$ on the right is a $K$-algebra (we identify
$K$ in $A$ with the elements of the form $(k,k)$ with $k\in K$) which is a regular ring,
and $$B(K^{[B_{1}]}\times K^{[B_{2}]})\cong B(K^{[B_{1}]})\times
B(K^{[B_{2}]}) \cong B_{1}\times B_{2}.$$ Making use of
Proposition~\ref{unicidad de algebras regulares},
we check that $A=conv(K)$. We know that $$A= K^{[B_{1}]}\times
K^{[B_{2}]}=conv(K)\times conv(K) =conv(K\times K)$$ so it
suffices to see that any $(k,k')\in K\times K$ is in $conv(K)$,
and this is trivial since $(k,k') = (1,0)(k,k) + (0,1)(k',k')$.
\end{dem}

\begin{lema}\label{iteracion}
    Let $A$ be a CFG-ring and $f:A\To A$ a polynomial map.
    Then, the set of iterated maps $\{f^{k} : k\in\mathbf{N}\}$ is finite.
\end{lema}

\begin{dem}
    First, it is easy to check that a polynomial map commutes with convex combinations:
    the identity map and constant maps do and this property is preserved under sums and products. 
    Suppose that $A=conv\{x_{1},\ldots,x_{n}\}$. For $i=1,\ldots,n$ we
    can express $f(x_{i})$ like a convex combination $f(x_{i}) =
    a_{1i}x_{1}+\cdots+a_{ni}x_{n}$. Call $R$ the subring of $B(A)$
    generated by the $a_{ij}$'s.  The ring $R$ is finite since it is
    a finitely generated subring of a Boolean ring. An induction argument shows that for all
    $k\in\mathbf{N}$, $f^{k}(x_{i})$ can be expressed as a convex
    combination with coefficients in $R$: if $f^{k-1}(x_{i})
    = b_{1}x_{1}+\cdots+b_{n}x_{n}$, then
    \begin{eqnarray*} f^{k}(x_{i}) &=& f(f^{k-1}(x_{i})) = f(\sum_{r=1}^{n}b_{r}x_{r})
    = \sum_{r=1}^{n}b_{r}f(x_{r})\\ &=&
    \sum_{r=1}^{n}b_{r}\sum_{t=1}^{n}a_{tr}x_{t} =
    \sum_{t=1}^{n}\left(\sum_{r=1}^{n}b_{r}a_{tr}\right)x_{t}.\end{eqnarray*}
    Therefore, for each $k\in\mathbf{N}$ there is a
    $n\times n$-matrix $M^{k}=\{a_{ij}^{k}\}$ of elements of $R$ such that
    $f^{k}(x_{j}) = \sum_{i}a_{ij}^{k}x_{i}$. Since $R$ is finite,
    the set of $n\times n$-matrices over $R$ is finite too.
\end{dem}

\begin{lema}\label{subanillo de un CFG}
Let $A$ be a CFG-ring and $R$ a finitely generated subring of
$A$. Then $R$ is finite.
\end{lema}

\begin{dem}
We proceed by induction on $n$, the number of generators. For
$n=0$, $R$ is the prime ring of $A$, that is finite, by applying
Lemma~\ref{iteracion} to the map $x\mapsto x+1$. Suppose the
assertion of the lemma for $n$, and we will prove it for $n+1$.
Any subring generated by $n+1$ elements is of the form $T[x]$ with
$T$ generated by $n$ elements and hence, by the induction
hypothesis, finite. $T[x] = \{\sum_{i=1}^{k}a_{i}x^{i} : a_{i}\in
T \}$. Applying again Lemma~\ref{iteracion} to $a\mapsto ax$, we
deduce that $\{x^{k} : k\in\mathbf{N}\}$ is finite, so $T[x]$ is
finite.
\end{dem}

\begin{lema}\label{ideales primos de una K-algebra regular}
Let $K$ be a field, $B$ a Boolean ring, and let $p$ be a prime
ideal of $A=K^{[B]}$. Then $A/p$ is isomorphic to $K$.
\end{lema}

\begin{dem} Let us check that the composition
$f:K\hookrightarrow A\To A/p$ is an isomorphism. It is injective,
since $K$ is a field. It is surjective: for all $x\in A$, express
$x$ as a convex combination of elements of $K$,
$x=\sum_{i=1}^{n}a_{i}k_{i}$. When mapping to the quotient ring,
all idempotents $a_{i}$ map into idempotents of the domain $A/p$,
that is, 0 or 1. So $x+p\in K/p = Im(f)$.
\end{dem}

\begin{teo}\label{estructura de CFG-anillos}
    Let $A$ be a CFG-ring. Then, there are  Boolean rings $B_{1},\dots,B_{n}$
    and non-isomorphic finite fields $K_{1},\dots,K_{n}$ such that
    $$A \cong K_{1}^{[B_{1}]} \times \cdots \times K_{n}^{[B_{n}]}$$
    Furthermore, this decomposition is unique, up to isomorphism and order.
\end{teo}

\begin{dem}
    Existence of such a decomposition: By Lemma~\ref{producto}, it
    suffices to find a decomposition where there may be isomorphic finite fields.
     Suppose $A = conv(H)$ where
    $H= \{x_{1},\dots,x_{n}\}$. Call $T$
    the subring generated by $H$, that is finite, by
    Lemma~\ref{subanillo de un CFG}. Therefore, $T$ is a finite reduced ring, so it is a
    product of finite fields (just apply
    the Chinese Remainder Theorem to the set of prime ideals of $T$). Take an isomorphism
    $h:K_{1}\times\cdots\times K_{m}\To T\hookrightarrow A$. Call
    $e_{i} = (\delta_{ij})_{j=1}^{m}$ ($\delta_ij$ is the Kronecker delta) and
    $\varepsilon_{i}=h(e_{i})$. We have a ring decomposition
    $A\cong \prod_{i=1}^{m}A\varepsilon_{i}$ since the
    $\varepsilon_{i}$'s constitute a c.f.o.i. of $A$.
     The restriction $h:K_{i}=Ke_{i}\To
    A\varepsilon_{i}$ provides a ring homomorphism, so we can view
    $A\varepsilon_{i}$ as a $K_{i}$-algebra, that is regular,
    since it is a factor of a regular ring. If we prove that
    $conv(h(K_{i})) = A\varepsilon_{i}=:A_{i}$ we will deduce, by
    Proposition~\ref{unicidad de algebras regulares}, that
    $A_{i} \cong K_{i}^{[B(A_{i})]}$. We have to
    see that any element of $A_{i}$ is a convex combination of
    elements of $h(K_{i})$ with scalars in $B(A_{i})$. Take $x\in
    A_{i}$. There is a convex combination in $A$,
    $x=\sum_{j}b_{j}r_{j}$ with $r_{j}=h(k_{j})\in T$. Then,
    $x=\varepsilon_{i}^{2}x
    =\sum_{j}(\varepsilon_{i}b_{j})(\varepsilon_{i}r_{j}) =
    \sum_{j}(\varepsilon_{i}b_{j})h(e_{i}k_{j})$ provides
    us the expression desired.\\

     Uniqueness: Suppose given one such
    a decomposition $A\cong A_{1}\times\cdots\times A_{n}$. Each
    $A_{i}=K_{i}^{[B_{i}]}$ can be seen as a principal ideal of $A$. Any prime ideal
    of $A$ is of the form $p_{i}^{e} = A_{1}\times\cdots\times
    A_{i-1}\times p_{i}\times A_{i+1}\times\cdots\times A_{n}$ for
    some prime ideal $p_{i}$ of $A_{i}$, and by Lemma~\ref{ideales
    primos de una K-algebra regular},
    $K_{i}\cong A_{i}/p_{i}\cong A/p_{i}^{e}$.
    Therefore, the fields $K_{i}$ are uniquely determined, up to
    isomorphism, by $A$, since they are those that appear as
    quotients by prime ideals. Furthermore, since $A_{i}$ is
    regular, the intersection of all
    prime ideals of $A_{i}$ is 0. Hence, the intersection of
    all prime ideals of $A$ whose quotients are isomorphic to
    $K_{i}$ is $A_{1}\times\cdots\times
    A_{i-1}\times 0\times A_{i+1}\times\cdots\times A_{n}$. Then,
    the intersection of all prime ideals of $A$ whose quotients are not
    isomorphic to $K_{j}$ is $0\times\cdots\times 0\times
    A_{j}\times\cdots\times 0\equiv A_{j}$. Therefore, the factor
    ring of the decomposition corresponding to $K_{j}$ is also
    uniquely determined by $A$, and also the Boolean ring $B_{j}$
    because it must be (isomorphic to) the ring of idempotents of
    that factor ring.
\end{dem}

\section{Other finiteness conditions}

During the proof of Lemma~\ref{iteracion}, it was observed that any 
polynomial map $f:A\To A$ commutes with convex combinations. Maps with 
this property and their relation with polynomials are object of study in~\cite{Avi}. 
We summarize in the following proposition the information that we need about this: 

\begin{prop}
Let $A$ be a regular ring. For a function $f:A\To A$ the following are 
equivalent:\begin{enumerate}
\item The map $f$ commutes with convex combinations, that is, for any 
$x_{1},\ldots x_{n}\in A$ and any c.f.o.i of $A$, $a_{1},\ldots,a_{n}$, 
$f(\sum_{1}^{n}a_{i}x_{i}) = \sum_{1}^{n}a_{i}f(x_{i})$.
\item $e(f(x)-f(y))\leq e(x-y)$ for all $x,y\in A$.
\end{enumerate}
Furthermore, if $A$ is a CFG-ring, then the following condition is also equivalent:
\begin{enumerate}
\setcounter{enumi}{2}
\item $f$ is polynomial map.\\
\end{enumerate} 
\end{prop}

A map verifying conditions $(1)$ and/or $(2)$ above will be called a contractive map. We observe
that the map $e:A\To A$ is always contractive. The equivalence of $(3)$ with the others is the
statement of~\cite[Lema 3.3]{Avi}.\\

We recall the four finiteness conditions exposed in the introduction:

(1) All quotients by prime ideals are finite.

(2) The cardinality of the quotients by prime ideals is bounded by an 
integer.

(2.5) All contractive maps $f:A\To A$ are polynomial maps.

(3) $A$ is a CFG-ring.\\

We complete the diagram of implications, except for $(2.5)\Rightarrow (3)$, 
which we do not know whether it is true or not. We begin with a characterization of
condition (2) which shows that (2.5) implies (2):\\
 
\renewcommand{\theenumi}{\roman{enumi}}
\begin{prop}\label{e contractiva}
For a regular ring $A$, the following are equivalent:
\begin{enumerate}
\item The map $e:A\To A$ is a polynomial map.
\item There exists a natural number $n$ such that $|A/p|<n$ for all prime
ideals $p$ of $A$.
\end{enumerate}
\end{prop}

\begin{dem}
Suppose that $e(x) = a_{k}x^{k}+\cdots +
a_{0}$. Take a prime ideal $p$ of $A$ and the natural projection
$\pi:A\To A/p$. Since $e(x)$ is idempotent, $\pi(e(x))$ is an
idempotent of the field $A/p$, so $\pi(e(x))\in\{0,1\}$ and

\begin{eqnarray*}0 &=& (\pi(e(x))-1)\pi(e(x))\\ &=&
(\pi(a_{k})\pi(x)^{k}+\cdots+\pi(a_{0})-1)(\pi(a_{k})\pi(x)^{k}+\cdots+\pi(a_{0}))\end{eqnarray*}
 for all $x\in A$. Therefore, all elements in the field $A/p$ are roots of a
polynomial of degree $2k$ and $|A/p|<2k$.\\

Conversely, suppose that $\{n_{1},\ldots,n_{r}\}$ is the set of
cardinalities of quotients of $A$ by prime ideals and call
$m_{i}=n_{i}-1$, $m=m_{1}\cdots m_{r}$, $k_{i}=m/m_{i}$. We prove that $e(x) = x^{m}$.
Since the intersection of all prime ideals is zero, it is sufficient to
check that $\pi(e(x)) = \pi(x)^{m}$ for every projection $\pi:A\To A/p$.
If $|A/p| = n_{i}$, then $$\pi(x)^{m} = \pi(x)^{m_{i}k_{i}} = (\pi(x)^{n_{i}-1})^{k_{i}} =
(1-\delta_{\pi(x)0})^{k_{i}} = 1-\delta_{\pi(x)0}$$
 where $\delta_{\pi(x)0}$ is the Kronecker delta. On the other hand, since
 $e(x)$ is an idempotent associated to $x$, $\pi(e(x))$ is an idempotent
 associated to $\pi(x)$ in the field $A/p$, so $\pi(e(x)) =
 1-\delta_{\pi(x)0}$.\\
\end{dem}

A counterexample for $(1\not\Rightarrow 2)$: Let $B$ be a Boolean ring and
$$I_{1}\supset I_{2}\supset I_{3} \supset\cdots$$ a
descending chain of ideals of $B$ such that for any prime ideal
$p$ of $B$ there is some $n$ with $p\supset I_{n}$.
Fix also a prime number
$q\in\mathbf{N}$ and a tower of finite fields of characteristic
$q$, $$F_{1}\subset F_{2}\subset F_{3}\subset\cdots$$ and call
$|F_{i}| = q^{n_{i}}$, $F=\bigcup_{1}^{\infty}F_{i}$. From these data, we 
will construct a ring $A=F_{*}^{[I_{*}]}$ that verifies $(1)$, and that if all 
ideals $I_{i}$ are nonzero, it does not verify $(2)$.
 For instance, we can take $B=\{0,1\}^{\mathbf{N}}$ and $$I_{n} = \{x\in B :
\{j: x_{j}=1\}\ is\ finite,\ x_{i}=0 \ for\ i<n\}.$$ (If $p$ is a
prime ideal of $B$, either it contains $I_{1}$ or it is of the
form $p_{i}=\{x: x_{i}=0\}$).\\

We define $A$ to be
the subset of $F^{[B]}$ formed by the convex combinations of
elements of $F$ for which the coefficients of elements that are
not in $F_{i}$ are in $I_{i}$.\\

 First we check that $A$ is a subring. It clearly contains $0$,
$1$ and $-1$. Suppose $x,y\in A$. We express $x$ and $y$ as convex
combinations of elements of $F$ like $x=\sum_{1}^{n}a_{i}x_{i}$
and $y=\sum_{1}^{m}b_{j}y_{j}$ in such a way that if $x_{i}\not\in
F_{k}$ then $a_{i}\in I_{k}$ and the same for $y$. Then
$x+y = \sum_{i,j}a_{i}b_{j}(x_{i}+y_{j})$. Now, if
$x_{i}+y_{j}\not\in F_{k}$ then either $x_{i}\not\in F_{k}$ or
$y_{j}\not\in F_{k}$. In any case, $a_{i}b_{j}\in I_{k}$. So
$x+y\in A$. The same reasoning leads to $xy\in A$.\\

 Since $F^{[B]}$ is
regular, $A$ is reduced. Also, $$B\cong
B(F^{[B]}) = conv\{0,1\}\subseteq A.$$

 We will see that $|A/p|$ is finite for all prime ideals of $A$.
That will prove also that $A$ is regular, by~\cite[Theorem 1.16]{Goodearl}.
Take $p$ a prime ideal of $A$. Then, $p\cap
B(A)$ is a prime ideal of $B(A)\cong B$, so $p\cap B(A)\supset
I_{k}$ for some $k$. Take $x\in A$ and express it as a convex
combination of elements of $F$, $x=\sum_{1}^{n}a_{i}x_{i}$ where
$a_i\in I_{k}\subset p$ whenever $x_i\not\in F_k$.
Hence, the class of any $x$ modulo $p$ is a linear combination of
classes of elements of $F_{k}$ with coefficients idempotents of
the domain $A/p$ and therefore $|A/p|\leq |F_{k}|<\infty$.\\

Finally, supposing all $I_{i}$ are nonzero, we will prove that for any $k>0$ there is a prime ideal
$p$ of $A$ such that $|F_{k}|\leq |A/p|$. Since $A$ is reduced,
the intersection of all prime ideals of $A$ is zero and there
exists a prime ideal $p$ of $A$ such that $0\neq
I_{k}A\not\subseteq p$. Take $a\in I_{k}$, $a\not\in p$. We claim
that the elements $\{ax : x\in F_{k}\}$ are in $A$ and are all
distinct modulo $p$. First, $ax =
ax+(1-a)0$ and provided $x\in F_{k}$, if $x\not\in F_i$ then $i>k$ and $a\in I_k\subset I_i$.
This shows the $ax\in A$ whenever $x\in F_k$. Now,
suppose $ax$ equals $ay$ modulo $p$ for $x,y\in F_{k}$, $x\neq y$.
Then $a(x-y)\in p$, and since $(x-y)^{-1}\in F_{k}$,
$a(x-y)^{-1}\in A$ and on the other hand, $p$ is an ideal in $A$,
so $a(x-y)a(x-y)^{-1} = a\in p$ and this leads to a
contradiction.\\

 A counterexample for $(2\not\Rightarrow 2.5)$: Let $K=\{0,1,a,b\}$ be a
field with four elements, and $A$ the subring of $K^{\mathbf{N}}$
formed by the sequences in which only a finite number of terms are
different from 0 and 1 (This example appears also in~\cite{Goodearl} for other 
purposes).\\

It is plain that $A$ is a reduced ring.\\

Since $x(x+1)(x^{2}+x+1)=0$ for all $x\in K$, the same relation
holds for all $x\in A$. If $p$ is a prime ideal of $A$, then $A/p$
is a domain and this formula implies that $|A/p|\leq 4$. In
particular, every $A/p$ is a field, so $A$ is a regular ring,
by~\cite[Theorem 1.16]{Goodearl}.\\

Finally, we find a contractive map $f:A\To A$ that is not
polynomial. Consider $\hat{a}\in K_{4}^{\mathbf{N}}$ the sequence
constant equal to $a$, and $f:A\To A$ given by $f(x) =
x(x+1)(x+\hat{a})$. The map $f$ is contractive,
since it is the restriction of a polynomial in
$K_{4}^{\mathbf{N}}$. Let us see now that $f$ is not a polynomial
map in $A$. Suppose it is: $f(x) = \sum_{i=0}^{m}c^{(i)}x^{i}$ for
some $c^{(i)}\in A$. Just by definition of $A$, there exists an
index $k$ such that $c^{(i)}_{k}\in\{0,1\}$ for all $i$. Then, for
all $x\in A$, we have $x_{k}(x_{k}+1)(x_{k}+a) = f(x)_{k} =
\sum_{i=0}^{m}c^{(i)}_{k}x^{i}_{k}$ with
$c^{(i)}_{k}\in\{0,1\}=\mathbf{Z}_{2}$. This is a contradiction
because the map $h:K_{4}\To K_{4}$ given by $h(t)=t(t+1)(t+a)$
cannot be given by a polynomial with coefficients in
$\mathbf{Z}_{2}$ since $h(a)=0$ but $h(b)\neq 0$.\\

One natural question is up to which point the strong structure theorem we 
have obtained for CFG-rings can be generalized to larger classes of reduced rings 
with finite quotient fields. We finish by making a remark and posing a 
problem in that direction.

\begin{prop}
Let $A$ be a reduced ring with all quotient fields finite. Then,
$char(A)=q_{1}\cdots q_{n}\neq 0$ with $q_{1},\ldots,q_{n}$ different prime 
numbers and there is a canonical decomposition $A\cong 
A_{1}\times\cdots\times A_{n}$ with $char(A_{i}) = q_{i}$.
\end{prop}

\begin{dem}
For each prime number $q>0$, since $A$ is a regular ring, $V(q) 
= \{p\in Spec(A) : q\in p\}$ is a closed open set. Since all quotient 
fields are finite, and hence have prime characteristic, the set of all 
$V(q)$ is an open covering of $Spec(A)$. By compactness $Spec(A)=V(q_{1})\cup\cdots\cup V(q_{m})$. Therefore $q_{1}\cdots 
q_{m}$ is in all prime ideals of $A$ and $char(A)$ divides $q_{1}\cdots 
q_{m}$.\\

 Since $V(q)\cap V(q') = \emptyset$ for different prime numbers 
$q,q'$, also $e(q)\vee e(q')=1\in A$, and $(1-e(q))(1-e(q'))=0$. On the other 
hand, since $Spec(A)=V(q_{1})\cup\cdots\cup V(q_{m})$, $$e(q_{1})\cdots 
e(q_{n})=0\ \ \text{ and }\ \ \bigvee_{1}^{m}(1-e(q_{i}))=1.$$ All this gives
$$A=(1-e(q_{1}))A\oplus\cdots\oplus (1-e(q_{n}))A$$ what induces the desired 
decomposition.
\end{dem}

\begin{prob}
Are all reduced rings with finite quotient fields of prime characteristic isomorphic 
to rings of type $F_{*}^{[I_{*}]}$ as constructed in the counterexample for $(1\not\Rightarrow 2)$? 
\end{prob}

If the answer were positive, there would be a complete structure theorem, not only for CFG-rings 
but for reduced rings verifying condition (2).

\end{document}